\documentclass[12pt]{amsart}
\usepackage[margin=1in]{geometry}
\usepackage{times}
\usepackage{amsthm,amssymb}
\usepackage{amsmath,xypic}
\usepackage[usenames,dvipsnames]{color}
\usepackage{natbib}
\let\cite\citep
\usepackage[colorlinks=true,hyperindex, citecolor=blue, urlcolor=cyan]{hyperref}
\usepackage{epigraph}
\usepackage{graphicx}
\usepackage{fancyhdr} \usepackage{enumitem}
\usepackage{footnote}
\usepackage{cleveref}
\newcommand{\be}{\begin{equation}}
\newcommand{\ee}{\end{equation}}
\newcommand{\bes}{\begin{equation*}}
\newcommand{\ees}{\end{equation*}}
\newcommand{\bea}{\begin{eqnarray}}
\newcommand{\eea}{\end{eqnarray}}
\newcommand{\beas}{\begin{eqnarray}}
\newcommand{\eeas}{\end{eqnarray}}
\newcommand{\ben}{\begin{note}}
\newcommand{\een}{\end{note}}
\newcommand{\bexl}{\vskip0.1em\noindent\hrulefill\vskip1em\begin{ExerciseList}}
\newcommand{\eexl}{\end{ExerciseList}\hrulefill}

\newcommand{\bthm}{\begin{theorem}}
\newcommand{\ethm}{\end{theorem}}
\newcommand{\bpro}{\begin{prop}}
\newcommand{\epro}{\end{prop}}
\newcommand{\bcor}{\begin{corollary}}
\newcommand{\ecor}{\end{corollary}}
\newcommand{\bcon}{\begin{conjecture}}
\newcommand{\econ}{\end{conjecture}}
\newcommand{\bp}{\begin{proof}}
\newcommand{\ep}{\end{proof}}
\newcommand{\blem}{\begin{lemma}}
\newcommand{\elem}{\end{lemma}}
\newcommand{\bn}{\begin{note}}
\newcommand{\en}{\end{note}}
\newcommand{\benum}{\begin{enumerate}}
\newcommand{\eenum}{\end{enumerate}}
\newcommand{\bed}{\begin{defn}}
\newcommand{\eed}{\end{defn}}
\newcommand{\brem}{\begin{remark}}
\newcommand{\erem}{\end{remark}}

\newcommand{\btik}{\begin{tikzpicture}\begin{axis}[scale=0.5,axis y line=center, axis x line=middle]}
\newcommand{\etik}{\end{axis}\end{tikzpicture}}

\let\into=\hookrightarrow
\let\mapsto=\longmapsto

\newcommand{\upperRomannumeral}[1]{\uppercase\expandafter{\romannumeral#1}}

\newtheorem{theorem}[equation]{Theorem}      \newtheorem{lemma}[equation]{Lemma}          \newtheorem{corollary}[equation]{Corollary}  \newtheorem{proposition}[equation]{Proposition}

\theoremstyle{definition}

\theoremstyle{definition}
\newtheorem{defn}[equation]{Definition}
\theoremstyle{remark}

\theoremstyle{definition}
\newtheorem{remark}[equation]{Remark}

\numberwithin{equation}{section}

\let\into=\hookrightarrow
\let\isom=\simeq
\let\rk=\rank
\let\tensor=\otimes

\newcommand{\Z}{{\mathbb Z}}

\renewcommand{\int}{\operatorname{int}}
\renewcommand{\O}{{\mathcal O}}
\renewcommand{\P}{{\mathbb P}}

\newcommand{\mapright}[1]{{\xymatrix{{}\ar[r]^{#1}&{}}}}

\renewcommand{\bpro}{\begin{proposition}}
\renewcommand{\epro}{\end{proposition}}

\newcommand{\pn}[1]{\mathbb{P}^{#1}}
\newcommand{\pnn}{\pn{n}}

\begin{document}

\title[]{On determinantal equations for curves and Frobenius split hypersurfaces}\author{Kirti Joshi}\address{Math. department, University of Arizona, 617 N Santa Rita, Tucson
85721-0089, USA.}

\thanks{}\subjclass{}\keywords{}

\begin{abstract}
I consider the problem of existence of intrinsic determinantal equations for plane projective curves and hypersurfaces in projective space in characteristic $p>0$ and prove that 
(1) if $X$ is an ordinary plane projective curve then there exists an intrinsic matrix $M$ with homogeneous linear entries such that $X$ is defined set theoretically by $\det(M)=0$. Additionally, if $p=2$, then $M$ is symmetric, and if $p=3$ then $M$ is skew-symmetric. (2) In higher dimensions and in any positive characteristic, any Frobenius split hypersurface in $\pnn$ is given by set theoretically as  the determinant of an intrinsic matrix with homogeneous entries of degrees between $1$ and $n-1$. In particular, this implies that any smooth, Fano hypersurface is set theoretically given by an intrinsic  determinantal equation and the same is also true for any Frobenius split Calabi-Yau hypersurface. In particular, in positive characteristic, a general smooth hypersurface of degree $\leq n+1$ in $\P^n$ is always given by an intrinsic set-theoretic determinantal equation.
\end{abstract}
\maketitle
\setlength{\epigraphwidth}{4in}
\begin{savenotes}
\epigraph{Ry\={o}kan! How nice to be like a fool\\ \  \ for then one's Way is grand beyond measure
}{(Master) Tainin Kokusen (to Ry\=okan) \citeauthor{ryokan}}
\end{savenotes}
\tableofcontents

\section{Introduction} The problem of finding determinantal equations for varieties goes back to the grandmasters of our subject. Readers should consult the two surveys (\cite{beauville00,beauville2018}) for an excellent introduction to this beautiful but difficult subject.

The classical problem  alluded to here  is the following. Let $k$ be an algebraically closed field. Let $X\into \P^n$ be a smooth projective hypersurface of degree $d$ given by an equation $G=0$. One says that $X$ is a \emph{determinantal hypersurface} if there exists a square matrix $M$,  whose entries are homogeneous linear polynomials in coordinates of $\pnn$ such that  $\det(M)=G$.   Note that this condition is invariant under automorphisms of $\pnn$. 

\newcommand{\sE}{\mathcal{E}}

More generally, one can ask, when does there exist a matrix $M$ of homogeneous polynomials in the coordinates of $\pnn$ and of degrees strictly less than that of $G$ such that $\det(M)=G^m$ for some integer $m\geq 1$.  If this weaker condition holds, then I will say that $X$ is a \emph{set theoretic determinantal hypersurface}.
\newcommand{\pfaff}[1]{{\rm Pfaff}(#1)}
If the entries of $M$  are linear (resp. quadratic, cubic, etc.) homogeneous polynomials then one says that $X$ is a linear (resp. quadratic, cubic, etc.) (set theoretic) determinantal hypersurface. If the matrix $M$ has entries of bounded degrees then one says that \emph{$X$ is a bounded (set theoretic) determinantal hypersurface}.  If the matrix $M$ is skew-symmetric and if $\pfaff M=G^m$ then one says that $X$ is a \emph{(set theoretic) Pfaffian determinantal hypersurface} (recall that the Pfaffian $\pfaff{M}$ of a skew-symmetric matrix satisfies $\det(M)=\pfaff{M}^2$).

One learns from (\cite{beauville00,beauville2018}) that, in modern parlance, the question of whether $X$ is or is not a determinantal hypersurface reduces to finding a coherent sheaf $\sE$ on the ambient projective space with a certain type of a (minimal) resolution;  the matrix $M=M_\sE$ is a part of this resolution datum. More precisely, the question reduces to finding an arithmetically Cohen-Macaulay coherent sheaf $\sE$ (with some additional properties) on $\P^n$. The matrix $M_\sE$ is well-defined up to a choice of coordinates and generators and, up to these choices, $M_{\sE}$ is an invariant of the isomorphism class of $\sE$. However, the isomorphism class of the  coherent sheaf $\sE$ which provides determinantal equations is seldom unique so there are many such (inequivalent) determinantal representations of $X$.

As is explained in (\cite[Proposition 1]{beauville00}) the problem of finding linear determinantal (resp. set theoretic linear determinantal) equations for $X$ is equivalent to finding an \emph{Ulrich line bundle} on $X$ (resp. \textit{Ulrich vector bundles} on $X$). 

As was proved in (\cite{eisenbud03}),  the existence of an Ulrich bundle on  any projective variety $X$ implies, remarkably, that the Chow form of $X$ is given (set-theoretically) by a single  linear determinantal equation in the Plucker coordinates on a suitable Grassmannian. The existence of Ulrich bundles is known in only in a few cases including complete intersections (so the conjecture of \cite{eisenbud03} is a wide open problem). On the other hand, the existence of Ulrich line bundles is  more restrictive even for hypersurfaces \cite[Proposition 7.6]{beauville00}, \cite[Section 3]{beauville2018}. 
 
This note concerns itself with the \emph{question of provenance} of  determinantal equations (resp. theoretic determinantal equations) when they exist: when does there exist an \emph{intrinsic} (set theoretic) determinantal equation for a hypersurface? Equivalently: when does there exist an \emph{intrinsic coherent sheaf} $\sE$ on $X$ which provides a (set theoretic) determinantal equation for $X$? 

I do not attempt to make a precise definition of an intrinsic sheaf because such an exercise might be too restrictive, but twists of embedding independent sheaves such as the tangent bundle, the cotangent bundle, the sheaves of differentials and their tensor powers, symmetric powers and Schur functors applied to these sheaves should be considered to be intrinsic sheaves on $X$ (the only embedding dependence is through the specific twist).

In characteristic zero, I do not know how to find intrinsic equations for hypersurfaces  and given the classical nature of the subject it would surprising to be able to find something new in this subject, and the results of \cite{raychaudhury2024}, \cite{raychaudhury2023-tangent} suggests that the problem is restrictive in characteristic zero. But my observation of this paper (and also of (\cite{joshi2019-ulrich} and \cite{joshi21,joshi22-corrigendum}) is that in characteristic $p>0$ this problem does seem to offer  a natural solution for curves and hypersurfaces under \emph{arithmetic} assumptions on $X$. 

The (intrinsic) sheaf $B^1_X$  which appears in my construction  is the canonically defined sheaf  of locally exact differentials (defined in \Cref{se:ord-split}) and it is coherent and even locally free  only in characteristic $p>0$ (for any smooth variety $X$). As is well-known, this sheaf plays an important role in the algebraic geometry in positive characteristics--and some of its important properties can be found in \cite[2.2]{Illusie1979},  \cite[Th\'eor\`eme 4.1.1]{Raynaud1982},  \cite[Proposition 1.1]{Illusie1990}, \cite{joshi04c}. The new discovery of (\cite{joshi2019-ulrich}) and  strengthened by the results of this paper is that the sheaf $B^1_X$ appears to play a some unexpected (in the extant work on syzygies) role  in syzygies of curves and higher dimensional varieties  which clearly needs to be understood better.

In Theorem~\ref{th:curves}, I show that for any smooth, projective, ordinary, plane curve over a field of positive characteristic there exist an intrinsic linear determinantal equation; for $p=2$ one may take this equation to be the determinant of an intrinsic symmetric matrix and  for $p=3$ one may take it to be an intrinsic skew-symmetric matrix i.e. a Pfaffian equation. Moreover, in Theorem~\ref{th:chow-form-for-curves} and for all $p\ge2$,  I show that the Chow form of any smooth, projective, ordinary curve (not necessarily  a plane curve) is always given by an intrinsic linear determinantal equation.  

Now let me mention results obtained here in higher dimensions. Theorem~\ref{th:frob-split-equations} shows that  under the arithmetic assumption of \emph{Frobenius splitting}  any smooth, projective Frobenius split hypersurface is given by an intrinsic set theoretic determinantal equation. In particular, one sees that any smooth, projective hypersurface in $\P^n$ of degree $\leq n$ is always given by an intrinsic set theoretic determinantal equation and any Frobenius split hypersurface of degree $n+1$ in $\P^n$ i.e. a Frobenius-split Calabi-Yau hypersurface (equivalently  an ordinary Calabi-Yau hypersurface \cite[Proposition 3.1(ii)]{joshi03}).  

\Cref{cor:general-hyp} and Corollary~\ref{cor:general-hyp-calabi-yau} show that a general hypersurface of degree $\leq n+1$  is an (intrinsic) set-theoretic determinantal hypersurface; especially a general smooth Calabi-Yau hypersurface in $\P^n$ is given by such a determinantal equation. On the other hand,
there is no obvious classical analogue of the results proved here.

It is a pleasure to thank N. Mohan Kumar for comments and corrections.

\section{Frobenius splitting and ordinarity}\label{se:ord-split}
Suppose $k$ is an algebraically closed field of characteristic $p>0$. Let $X$ be a smooth, projective variety of dimension $n$. Let $F:X\to X$ be the absolute Frobenius morphism of $X$ (see \cite[Page 301]{Hartshorne1977}) which, on local sections $f$ of the structure sheaf $\O_X$, is given by $f\mapsto f^p$. 
Let $d:\O_X\to \Omega^1_X$ be the differential. One wants to consider the sheaf $B_X^1=d(\O_X)\subset \Omega^1_X$ \cite[4.1]{Raynaud1982}  of locally exact differentials. One of the features of working in characteristic $p>0$ and with a smooth variety $X$ is that this sheaf of locally exact differentials is locally free of finite rank. [This is never the case in characteristic zero.] The kernel of $d$ consists of local sections of $\O_X$ of the form $f^p$. Moreover, for local sections $f,g$ of $\O_X$, one has $d(f^pg)=f^pdg$ from which one sees that $B^1_X$ is a (locally free) subsheaf $B^1_X\subset F_*(\Omega^1_X)$ and  the differential $d:\O_X\to \Omega^1_X$ provides  the fundamental exact sequence
\be\label{eq:fundamental-seq} 
0\to \O_X\to F_*(\O_X)\to B^1_X\to 0. \ee
Note that if $X$ is a smooth projective curve then $B^1_X$ is locally free of rank $p-1$ and degree $(p-1)(g-1)$ \cite[4.1]{Raynaud1982}.  

A smooth, projective variety over $k$ is said to be \emph{Frobenius split} if and only if \eqref{eq:fundamental-seq} splits as an exact sequence of $\O_X$-modules. 
By (\cite[Proof of Prop. 3.1]{joshi03}) it is immediate from Frobenius splitting hypothesis that one has
\be\label{eq:vanishing}
H^i(X,B^1_X)=0 \textrm{ for all }i\geq 0.
\ee

One says that a smooth projective curve $X$ is  \emph{ordinary} if $H^i(X,B^1_X)=0$ for all $i\geq 0$ (equivalently Frobenius morphism $H^i(X,\O_X)\to H^i(X,\O_X)$ is an isomorphism for all $i\geq0$. [The notion of ordinarity extends to all dimensions \cite[D\'efinition 1.1]{Illusie1990} but is not needed here.]
\brem\ 
Note that
\benum[label={(\bf\arabic{*})}]
\item A smooth, projective curve $X$ is ordinary  if and only if  the Jacobian of $X$ is an ordinary abelian variety in the sense of \cite[Definition 1.1]{Illusie1990}. 
\item If a smooth, projective curve of genus $g$ is Frobenius split, then $g\leq 1$.
\item On the other hand, ordinary curves of all genera exist \cite{Illusie1990}.
\item For a comparison of Frobenius splitting and ordinarity see \cite{joshi03}.
\eenum
\erem

The following lemma is standard and is useful in understanding the Frobenius splitting hypothesis.
\blem\label{le:frob-deg} 
Let $X\subset\pnn$ be a smooth, projective and Frobenius split hypersurface. Then
$$\deg(X)\leq n+1.$$
\elem
\bp 
Let $\omega_X$ be the canonical line bundle of $X$. By \cite[Section 2]{mehta1985}, Frobenius splitting, if it exists, gives a section $0\neq s\in H^0(X,\omega_X^{1-p})$ and hence this cohomology group is non-zero. Since $\omega_X=\O_X(d-n-1)$, one sees that $H^0(X,\omega_X^{1-p})\neq 0$ if and only if $(1-p)(\deg(X)-n-1)\geq0$ and since $p\geq 2$, the assertion of the lemma follows.
\ep

\brem\label{re:open-ness}\ 
\benum[label={(\bf\arabic{*})}]
\item If $X\to S$ is a flat family of smooth, proper varieties, then there is an open (possibly empty)  subset   of $S$ over which geometric fibers are Frobenius split (this is immediate from the \cite[Proposition 9]{mehta1985}).
\item Similarly, if $X\to S$ is a smooth, proper family of curves, then there is an open subset of $S$ over which the fibers are ordinary \cite[Proposition 1.2]{Illusie1990}. 
\item If $S$ is, say, the moduli of space curves, or scheme parameterizing smooth hypersurfaces in $\P^n$ of degree at most $n+1$, then this open subset of $S$ given by (1,2) is also non-empty (see \cite[Th\'eor\`eme 2.2]{Illusie1990}). 
\item Hence, Frobenius splitting (resp. ordinarity) can be viewed as an \emph{arithmetic condition} as well as  a \emph{geometric genericity condition} in moduli. 
\item However, it is important to note that Frobenius splitting (resp. ordinarity) is a condition which can often be tested for a given curve. For example, the Fermat curve $x^n+y^n+z^n=0$ is ordinary if and only if $p\equiv1 \bmod{n}$ \cite[2.3.3]{Illusie1990}, whereas most genericity conditions are often difficult to check for a given curve or a hypersurface. 
\eenum
\erem

\section{Determinantal equations}

\newcommand{\sF}{\mathcal{F}}
\newcommand{\sG}{\mathcal{G}}
\newcommand{\rk}{{\rm rk}}
\newcommand{\reg}{{\rm reg}}

A vector bundle $\sE$ on $\pnn$ will called a \emph{lineal bundle} if $\sE$ is a direct sum of line bundles on $\pnn$.  Thus, $\sE$ is a lineal bundle on $\pnn$ if and only if $\sE\isom \bigoplus_i\O_{\P^n}(a_i)$ for some $a_i\in \Z$.

Let $X\into\P^n$ be a smooth, projective hypersurface of dimension $\geq 1$ (so $n\geq 2$). The following result combines (\cite[Theorem A]{beauville00} and \cite[Proposition 1]{beauville2018}) and establishes a correspondence between sheaves of appropriate sort on $\pnn$ and the existence of set theoretic  (resp. linear) determinantal equations for $X$. 
\bpro\label{p:criterion}
Let $X\into \pnn$ be a smooth, projective hypersurface given by an equation $G=0$ of degree $d$ and suppose $\sE$ is a coherent sheaf on $\pnn$ supported  on $X$.
\benum[label={(\bf\arabic{*})}]
\item  The  following are equivalent:
\benum[label={(\roman{*})}]
\item $\sE$ is an ACM bundle on $X$ (i.e. $\sE$ is a vector bundle on $X$ which an ACM as an $\O_{\pnn}$-module).
\item $\sE$ admits a minimal resolution of the form
$$0 \to \sF_1\mapright{M} \sF_0 \to \sE \to 0,$$
where $\sF_1,\sF_0$ are  lineal bundles of rank $\rk(\sF_1)=r=\rk(\sF_0)$ on $\pnn$.
\item  $\det(M)=G^{r}$, in other words, $X$ is a set theoretic determinantal hypersurface.
\eenum
\item The following are equivalent:
\benum[label={(\roman{*})}]
\item $\sE$ is an Ulrich bundle on $X$.
\item $\sE$ admits a minimal linear resolution of the form
$$0 \to \sF_1=\O_{\P^n}(-1)^{rd}\mapright{M} \sF_0 =\O_{\P^n}^{rd} \to \sE \to 0,$$
where the entries of $M$ are homogeneous linear forms in the coordinates of $\P^n$,
\item  $\det(M)=G^{r}$, and $M$ has homogeneous linear entries i.e. $X$ is a set theoretic linear determinantal hypersurface.
\eenum
\eenum
\epro
\bp 
The assertion {\bf(1)} is \cite[Theorem A]{beauville00} and {\bf(2)} is \cite[Proposition 1]{beauville2018}.
\ep

This proposition reduces the task of finding set theoretic determinantal equations to finding vector bundles with appropriate kind of resolutions.

\section{Intrinsic determinantal equations for plain curves}
\bthm\label{th:curves}
Let $X$ be a smooth, projective, ordinary curve $X\subset\P^2$ over an  algebraically closed field $k$ of characteristic $p>0$. 
\benum[label={\bf{(\arabic{*})}}]
\item\label{th:curves-p5} Suppose $p\geq2$ then $X$ is an intrinsically  set theoretic linear determinantal curve.
\item\label{th:curves-p2} Suppose $p=2$ then $X$ is an intrinsically linear, symmetric determinantal curve.
\item\label{th:curves-p3} Suppose $p=3$ then $X$ is an intrinsically   linear Pfaffian determinantal curve.
\eenum
\ethm
\bp 
Since $X$ is ordinary,  by \cite{joshi2019-ulrich} one sees that $B^1_X(1)$ is an Ulrich bundle on $X$ and ordinarity of $X$ also provides the vanishing 
\be\label{eq:vanishing2} 
H^i(X,B^1_X)=0 \text{ for all } i\geq 0.
\ee
Note that \ref{th:curves-p5} is immediate from \cite{joshi2019-ulrich} which says that if $X$ is ordinary, then $B^1_X(1)$ is an Ulrich bundle on $X$ (and hence \Cref{p:criterion}{\bf(2)} gives the result). Since $B^1_X$ is an intrinsically defined vector bundle on $X$, one sees that $X$ is intrinsically a set theoretic linear determinantal curve in $\P^2$.

So it remains to prove the remaining two assertions. Let $\omega_X=\Omega^1_X$ be the canonical divisor of $X$. From \cite[4.1]{Raynaud1982} one sees that  $B^1_X$ carries a natural, perfect pairing  
\be\label{eq:pairing} B^1_X\tensor B^1_X\to \omega_X\ee which is skew-symmetric if $p$ is odd and symmetric if $p=2$.

Thus to prove \ref{th:curves-p2} it suffices to observe that for $p=2$, $B^1_X$ is a line bundle of degree $\deg(B^1_X)=g-1$ and the pairing \eqref{eq:pairing} gives  $(B^1_X)^2=\omega_X$, so $B^1_X$ is a natural theta divisor on $X$ \cite[4.1]{Raynaud1982}. Ordinarity of $X$ gives $H^0(X,B^1_X)=0$ and hence \ref{th:curves-p2} now follows from \cite[Proposition~4.2]{beauville00}.

Now suppose $p=3$, then  ordinarity of $X$ gives $H^0(X,B^1_X)=0$, $B^1_X$  is a locally free sheaf of rank two equipped with a natural skew-symmetric pairing \eqref{eq:pairing}. Hence one has a canonical isomorphism $$\det(B^1_X)=\wedge^2 B^1_X=\omega_X.$$   Now assertion \ref{th:curves-p3} follows from \cite[Proposition 5.1]{beauville00}. 

This completes the proof of \Cref{th:curves}.
\ep

\brem 
Let $p=3$. Let $m\geq 4$ be an integer. Then the Fermat curve given by the polynomial $G=0$ where $G=x^m+y^m+z^m$ in $\P^2$ is ordinary if and only if $p\equiv 1\bmod{m}$ \cite[2.3.3]{Illusie1990}. Thus \Cref{th:curves}\ref{th:curves-p3}  asserts that, for some integer $\ell\geq1$, there exists an intrinsic (in the above sense) $2\ell\times2\ell$ skew-symmetric  matrix $M$  with homogeneous linear entries such that $\det(M)=G^2$ and $G=\text{Pfaff}(M)$. From these properties, or using \cite[Proposition 5.1]{beauville00}, one obtains $\ell=m$. Thus $M$ is $2m\times 2m$ matrix with asserted properties. Using computer algebra software one can make $M$ explicit for small values of $m$ -- this exercise is left to the readers.
\erem

\section{Intrinsic Chow form for ordinary curves in projective space}
Suppose $X\into\P^n$ is a smooth, projective curve embedded in $\P^n$ and equipped with the ample line bundle $\O_X(1)$ given by this embedding. Let ${\rm Grass}(2,n)$ be the Grassmannian of lines in $\P^n$. Then the image of the incidence scheme $${\rm Inc}(X)=\left\{(x,\ell)\in X\times {\rm Grass}(2,n):x\in \ell\right\}$$ under the projection to ${\rm Grass}(2,n)$ is a divisor, denoted ${\rm Chow}(X)$, and called the \emph{Chow divisor} or the \emph{Chow form} of $X$ in $\P^n$. To put it differently, the Chow form of $X$ is the scheme of lines in $\P^n$ which meet $X$. More generally, as is shown in (\cite[Theorem~1.4]{eisenbud03}), if $\sE$ is a vector bundle on $X$ then $\sE$ has a Chow divisor, denoted ${\rm Chow}(\sE)$ in ${\rm Grass}(2,n)$ which satisfies $${\rm Chow}(\sE)=\rk(\sE)\cdot {\rm Chow}(X).$$ Combining this result with the (\cite[Theorem~2.1]{joshi2019-ulrich}) one has the following 
\bthm\label{th:chow-form-for-curves}
Let $X\into\P^n$ be a smooth, projective ordinary curve then $X$ has an intrinsic set-theoretic linear determinantal Chow form. More precisely there exists a divisor  $D\into {\rm Grass}(2,n)$ which is given globally by a single intrinsic linear determinantal equation in coordinates of ${\rm Grass}(2,n)$ and whose support is the Chow form of $X$.
\ethm

\newcommand{\sL}{L}
\newcommand{\boxten}{\mathop\boxtimes\displaylimits}

\section{Regularity of $B^1_X$ for Frobenius split varieties}
The main theorem of this section (Theorem~\ref{th:reg-bound}) provides a bound on the Mumford-Castelnuovo regularity of $B^1_X$  for any smooth, projective Frobenius split variety over an algebraically closed field in characteristic $p>0$. This result will be needed in the next section. Recall that if $\sF$ is a coherent sheaf on $\pnn$ then its Mumford-Castelnuovo regularity, denoted $\reg(\sF)$, is the smallest integer $d$ such that 
$$H^i(\pnn,\sF(d-i))=0\qquad\forall i>0.$$ 
\bthm\label{th:reg-bound}
Let $X\into \pnn$ be a smooth projective variety over an algebraically closed field characteristic $p>0$ equipped with $\O_X(1)$ provided by this embedding.
Assume that $X$ is Frobenius split.
Then $\reg(B^1_X)\leq \dim(X)$.
\ethm
\bp 
So one has to prove that $H^i(\pnn,B^1_X(d-i))=0$ for all $i>0$ and $d=\dim(X)$, or equivalently (by Leray spectral sequence for $X\into\pnn$) one has to prove that $H^i(X,B^1_X(d-i))=0$ for $i>0$. This is proved as follows. First note that if $i=d$ then $H^d(X,B^1_X)=0$  by \eqref{eq:vanishing} which is given by the assumption that $X$ is Frobenius split. 

So assume $1\leq i\leq d-1$. Twisting the split exact sequence \eqref{eq:fundamental-seq} (this is where Frobenius splitting is used again) by $\O_X(d-i)$ one has a split exact sequence
\be
0\to \O_X(d-i)\to F_*(\O_X)(d-i)\to B^1_X(d-i)\to 0,
\ee
where $1\leq i\leq d-1$. So $H^i(B^1_X(d-i))$ is a direct summand of $H^i(F_*(\O_X)(d-i))$ and so it is enough to prove that the latter  is zero
for $1\leq i\leq d-1$.
For $i$ in this range, $\O_X(d-i)$ is ample and $$H^i(X,F_*(\O_X)(d-i))=H^i(X,\O_X(p(d-i)))$$ by the projection formula. By (\cite[Proposition 1]{mehta1985}) as $X$ is Frobenius split one sees that $$H^i(X,\O_X(p(d-i)))=0$$ and hence its direct summand $H^i(X,B^1_X(d-i))=0$ for $1\leq i\leq d-1$. Thus the result is established.
\ep

\section{Determinantal Equations for Frobenius split Hypersurfaces}
Let me now prove the following theorem. 
\bthm\label{th:frob-split-equations}
Let $X\into\pnn$ be a smooth, projective, Frobenius split hypersurface  given by some equation $G=0$. Then there exists an intrinsic square matrix $M$ whose non-zero entries  have degrees bounded between  $1$ and $n-1$ and such that $\det(M)=G^r$ for some $r\geq 0$. In other words, every smooth, projective Frobenius split hypersurface is always  an intrinsic, bounded set theoretic determinantal hypersurface. In particular, every smooth Fano hypersurface is a bounded set theoretic determinantal hypersurface and every Frobenius split  Calabi-Yau hypersurface is a bounded set theoretic determinantal hypersurface.
\ethm

\brem\ 
\benum[label={(\bf\arabic{*})}]
\item I do not know any classical analogue of this theorem.
\item By \Cref{le:frob-deg}, for any $X$ as in \Cref{th:frob-split-equations}, one has $\deg(X)\leq n+1$ i.e. $X$ is Fano or a Calabi-Yau hypersurface. 
\item By (\cite{Fedder1987}) one sees that any smooth Fano hypersurface in $\pnn$ is Frobenius split. 
\item The Fermat quartic surface $x^4+y^4+z^4+w^4=0$ in $\P^3$ with $p\equiv3\bmod{4}$ is not Frobenius split \cite[2.3.3]{Illusie1990}.
\item For more examples of Frobenius split varieties see (\cite{joshi03}, \cite{joshi2019-ulrich}).
\eenum
\erem

The following Lemma will be used in the proof.
\renewcommand{\hom}{{\rm Hom}}
\blem\label{le:b1-hom}
Let $X\into \pnn$ be a smooth projective variety over an algebraically closed field characteristic $p>0$ equipped with $\O_X(1)$ provided by this embedding.
Assume that $X$ is Frobenius split with $\dim(X)\geq 2$. Then 
$$\hom(\O_X(m),B^1_X)=0\text{ for all }m\geq0.$$
\elem
\bp 
First observe that if $m=0$, then $\hom(\O_X,B^1_X)=0$ is equivalent to $H^0(X,B^1_X)=0$ which is immediate from \eqref{eq:fundamental-seq} as $X$ is Frobenius split. Next assume $m\geq 1$ then (using the fact that $X$ is Frobenius split) $$\hom(\O_X(m),B^1_X)=H^0(X,B^1_X(-m))\subseteq H^0(X,F_*(\O_X)(-m))=H^0(\O_X(-pm))=0$$ as $\O_X(-pm)$ is anti-ample for $m\geq1$. This proves the lemma. 
\ep

\bp[Proof of \Cref{th:frob-split-equations}] 
The regularity bound \Cref{th:reg-bound} is the  reason why one gets bounded degrees in \Cref{th:frob-split-equations}. Let $X\subset\pnn$ be a smooth, projective hypersurface. Then by Theorem~\ref{th:reg-bound} one has $\reg(B^1_X)\leq \dim(X)= n-1$. Let 
\be \label{eq:res} 0 \to \sG_1=\oplus_i\O_{\pnn}(-b_i)\mapright{M} \sG_0 \to B^1_X \to 0\ee
be a minimal resolution of $B^1_X$. By (\cite[Theorem 3.1]{joshi2019-ulrich}) $B^1_X$ is arithmetically Cohen-Macaulay, coherent sheaf on $\pnn$ one can read off the regularity of $B^1_X$ from the $b_i$. Specifically one has (see \cite[Exercise 4, Page 85]{eisenbud-syzygies-book})
$$\reg(B^1_X)=\max(b_i)-1.$$
As $\reg(B^1_X)\leq n-1$, by Theorem~\ref{th:reg-bound} one sees that $\max(b_i)-1\leq n-1$ (by Theorem~\ref{th:reg-bound}). Hence $\max(b_i)\leq n-1+1=n$. Hence $b_i\leq n$ for all $i$. 

I claim that $\hom(\O_{\P^n}(m),\sG_0)=0$ for any $m\geq 0$. 
Indeed if this is non-zero then one has a non-zero morphism $\O_X(m)\to B^1_X$ with $m\geq 0$ which is impossible by Lemma~\ref{le:b1-hom}. Since $\sG_0$ is a direct sum of line bundles of the form $\O_{\pnn}(a_j)$, one sees from this that $a_j\leq -1$. One sees from this that $H^0(\pnn,\sG_0)=0$. Moreover, $\hom(\O_{\P^n}(m),\sG_0)=0$ says that $\O_X(m)\subseteq \ker(\sG_0\to B^1_X)=\sG_1$ and so one has a trivial subcomplex $0\to\O_X(m)=\O_X(m)\to 0$ of our minimal resolution which contradicts minimality of the resolution. Alternately, one sees from \eqref{eq:res} that $0\to H^0(\pnn,\sG_1)\into H^0(\pnn,\sG_0)=0$.  Thus $\hom(\O_X(m),\sG_1)=0$ for any $m\geq 0$. Thus any line bundle which is a free direct summand of $\sG_0$ or $\sG_1$ is of the form $\O_{\pnn}(-m)$ with $m\geq 1$. 
Thus $B^1_X$ has a minimal resolution of the asserted type where the entries of the matrix $M$  for the morphism $\sG_1\to \sG_0$ are homogeneous polynomials in coordinates of $\P^n$ of   degrees bounded by $1$ and $n-1$.
\ep

\bcor\label{cor:general-hyp}
A general hypersurface of degree $\leq n+1$ in $\P^n$ ($n\geq 2$) is the support of an intrinsic determinantal equation arising from a matrix whose non-zero entries are homogenous  of degrees bounded by $1$ and $n-1$. 
\ecor
\bp
To prove that a general hypersurface of degree $\leq n+1$ in $\P^n$ is the support of an intrinsic determinantal equation, it is enough to observe that Frobenius splitting is an open condition in the space of hypersurfaces of a fixed degree \Cref{re:open-ness}{\bf(1)}. If degree $<n+1$ then by (\cite{Fedder1987}) then any such hypersurface is Frobenius split. If the degree is equal to $n+1$ then by (\cite[Th\'eor\`eme 2.2]{Illusie1990}) a general hypersurface is ordinary and by  (\cite{joshi03}) an ordinary Calabi-Yau hypersurface (i.e. of degree $n+1$) is Frobenius split. 
\ep

\bcor\label{cor:general-hyp-calabi-yau}
Any general Calabi-Yau hypersurface in $\P^n$ ($n\geq 2$) is the support of an intrinsic determinantal equation arising from a matrix whose non-zero homogenous entries have degrees bounded by $1$ and $n-1$. 
\ecor

\bibliographystyle{plainnat}
\bibliography{../../master/masterofallbibs.bib,
../../master/joshi.bib
}

\begin{thebibliography}{18}
\providecommand{\natexlab}[1]{#1}
\providecommand{\url}[1]{\texttt{#1}}
\expandafter\ifx\csname urlstyle\endcsname\relax
  \providecommand{\doi}[1]{doi: #1}\else
  \providecommand{\doi}{doi: \begingroup \urlstyle{rm}\Url}\fi

\bibitem[Ab\'e and Haskel(1996)]{ryokan}
Ryuichi Ab\'e and Peter Haskel.
\newblock \emph{The great fool: Zen master Ry\=okan Poems, Letters and other
  Writings}.
\newblock University of Hawaii Press, 1996.

\bibitem[Antonelli et~al.(2024)Antonelli, Casnati, Lopez, and
  Raychaudhury]{raychaudhury2024}
Vincenzo Antonelli, Gianfranco Casnati, Angelo~Felice Lopez, and Debaditya
  Raychaudhury.
\newblock On varieties with ulrich twisted conormal bundles, 2024.
\newblock URL \url{https://arxiv.org/abs/2306.00113}.

\bibitem[Beauville(2000)]{beauville00}
Arnaud Beauville.
\newblock Determinantal hypersurfaces.
\newblock \emph{Michigan Math. J.}, 48:\penalty0 39--64, 2000.
\newblock URL \url{https://doi.org/10.1307/mmj/1030132707}.

\bibitem[Beauville(2018)]{beauville2018}
Arnaud Beauville.
\newblock An introduction to ulrich bundles.
\newblock \emph{European Journal of Mathematics}, 4\penalty0 (1):\penalty0
  26--36, 2018.
\newblock URL \url{https://doi.org/10.1007/s40879-017-0154-4}.

\bibitem[Eisenbud(2005)]{eisenbud-syzygies-book}
David Eisenbud.
\newblock \emph{The geometry of syzygies}, volume 229 of \emph{Graduate Texts
  in Mathematics}.
\newblock Springer-Verlag, New York, 2005.
\newblock A second course in commutative algebra and algebraic geometry.

\bibitem[Eisenbud et~al.(2003)Eisenbud, Schreyer, and Weyman]{eisenbud03}
David Eisenbud, Frank-Olaf Schreyer, and Jerzy Weyman.
\newblock Resultants and {C}how forms via exterior syzygies.
\newblock \emph{J. Amer. Math. Soc.}, 16\penalty0 (3):\penalty0 537--579, 2003.

\bibitem[Fedder(1987)]{Fedder1987}
R.~Fedder.
\newblock {$F$}-purity and rational singularity in graded complete intersection
  rings.
\newblock \emph{Trans. {A}mer. {M}ath. {S}oc.}, 301\penalty0 (1):\penalty0
  47--62, 1987.

\bibitem[Hartshorne(1977)]{Hartshorne1977}
R.~Hartshorne.
\newblock \emph{Algebraic Geometry}.
\newblock Number~52 in Graduate {T}exts in {M}athematics. Springer-{V}erlag,
  New {Y}ork-{Heidelberg}, 1977.

\bibitem[Illusie(1990)]{Illusie1990}
L.~Illusie.
\newblock \emph{The {G}rothendieck {F}estschrift, Vol 2.}, volume~87 of
  \emph{Progr. {M}ath.}, chapter Ordinarit\'e des intersections compl\`etes
  g\'en\'erales, pages 376--405.
\newblock Bikhauser, 1990.

\bibitem[Illusie(1979)]{Illusie1979}
Luc Illusie.
\newblock Complexe de de {R}ham-{W}itt.
\newblock \emph{Ast\'erisque}, 63:\penalty0 83--112, 1979.

\bibitem[Joshi(2004)]{joshi04c}
Kirti Joshi.
\newblock Stability and locally exact differentials on a curve.
\newblock \emph{C. R. Math. Acad. Sci. Paris}, 338\penalty0 (11):\penalty0
  869--872, 2004.
\newblock URL \url{http://dx.doi.org/10.1016/j.crma.2004.02.019}.

\bibitem[Joshi(2019)]{joshi2019-ulrich}
Kirti Joshi.
\newblock A remark on ulrich and acm bundles.
\newblock \emph{Journal of Algebra}, 527:\penalty0 20--29, 2019.
\newblock ISSN 0021-8693.
\newblock \doi{https://doi.org/10.1016/j.jalgebra.2019.01.016}.
\newblock URL
  \url{https://www.sciencedirect.com/science/article/pii/S002186931930050X}.

\bibitem[Joshi(2021)]{joshi21}
Kirti Joshi.
\newblock On the construction of weakly ulrich bundles.
\newblock \emph{Advances in Mathematics}, 381:\penalty0 107598, 2021.
\newblock \doi{https://doi.org/10.1016/j.aim.2021.107598}.
\newblock URL
  \url{https://www.sciencedirect.com/science/article/pii/S0001870821000360}.

\bibitem[Joshi(2022)]{joshi22-corrigendum}
Kirti Joshi.
\newblock Corrigendum to “on the construction of weakly ulrich bundles”
  [adv. math. 381 (2021) 107598].
\newblock \emph{Advances in Mathematics}, 394:\penalty0 108025, 2022.
\newblock \doi{https://doi.org/10.1016/j.aim.2021.108025}.

\bibitem[Joshi and Rajan(2003)]{joshi03}
Kirti Joshi and C.~S. Rajan.
\newblock Frobenius splitting and ordinarity.
\newblock \emph{Int. Math. Res. Not.}, \penalty0 (2):\penalty0 109--121, 2003.
\newblock URL \url{http://dx.doi.org/10.1155/S1073792803112135}.

\bibitem[Lopez and Raychaudhury(2023)]{raychaudhury2023-tangent}
Angelo~Felice Lopez and Debaditya Raychaudhury.
\newblock On varieties with ulrich twisted tangent bundles, 2023.
\newblock URL \url{https://arxiv.org/abs/2301.03104}.

\bibitem[Mehta and Ramanathan(1985)]{mehta1985}
V.~B. Mehta and A.~Ramanathan.
\newblock {F}robenius splitting and cohomology vanishing of {S}chubert
  varieties.
\newblock \emph{Annals of Math.}, 122:\penalty0 27--40, 1985.

\bibitem[Raynaud(1982)]{Raynaud1982}
M.~Raynaud.
\newblock Sections des fibr\'es vectoriels sur une courbe.
\newblock \emph{Bull. {S}oc. {M}ath. {F}rance}, 110:\penalty0 103--125, 1982.

\end{thebibliography}
\end{document}